\documentstyle[12pt]{article}
\newcommand{\const}{\mathop{\rm const}\limits}

\newcommand{\grad}{\mathop{\rm grad}\limits}

\newcommand{\Rad}{\mathop{\rm Rad}\limits}

\begin{document}
\begin{center}

{\bf Exact constant in Sobolev's  and Sobolev's trace inequalities }\\

\vspace{3mm}

  {\bf for Grand Lebesgue Spaces }\\

\vspace{3mm}

{\sc Ostrovsky E., Rogover E. and Sirota L.}\\

\normalsize

\vspace{3mm}
{\it Department of Mathematics and Statistics, Bar-Ilan University,
59200, Ramat Gan, Israel.}\\
e-mail: \ galo@list.ru \\

\vspace{3mm}

{\it Department of Mathematics and Statistics, Bar-Ilan University,
59200, Ramat Gan, Israel.}\\
e - mail: \ rogovee@gmail.com \\

\vspace{3mm}

{\it Department of Mathematics and Statistics, Bar-Ilan University,
59200, Ramat Gan, Israel.}\\
e - mail: \ sirota@zahav.net.il \\

\end{center}

\vspace{4mm}

 {\it Abstract.} In this article we generalize the classical
Sobolev's  and Sobolev's trace inequalities  on the Grand
 Lebesgue Spaces instead the classical Lebesgue Spaces. \par
 We will distinguish the classical Sobolev's inequality and the so-called
 trace Sobolev's inequality. \par
  We consider for simplicity only the case of whole space.\par

 \vspace{3mm}

 {\it Key words:} Sobolev's and Poincare's inequalities, derivative, gradient, norm,
 Lebesgue spaces, Talenti's estimate, Bilateral Grand Lebesgue  spaces, trace,
 counterexamples. \par

 \vspace{3mm}

{\it Mathematics Subject Classification (2000):} primary 60G17; \ secondary
 60E07; 60G70.\\

\vspace{3mm}

\section{Introduction. Notations. Statement of problem.}

\vspace{3mm}
{\bf A. Ordinary Sobolev's inequality.} \par
 The classical Sobolev's inequality in the whole space $ R^m, $ see, e.g.
 \cite{Kantorovicz1}, chapter 11, section 5; \cite{Sobolev1}, \cite{Talenti1} etc.
 asserts that for all function $ f, f: R^m \to R, \ m \ge 3 $ from the Sobolev's space
 $ W^1_p(R^m), $ which may be defined as a closure in the Sobolev's norm

 $$
 ||f||W^1_p(R^m)= |f|_p + |Df|_p
 $$
  of the set of all finite continuous differentiable functions $ f, f: R^m \to R, $ that

  $$
  |f|_{q} \le K_m(p) \ |Df|_p, \ q = q(p)= mp/(m-p), \ p \in [1,m), \
  q \in (m/(m-1), \infty).\eqno(1)
  $$
  Here  $ m = 3,4, \ldots; $

  $$
  |f|_p = |f|_{p,m} = |f|_{p,R^m} = \left[ \int_{R^m} |f(x)|^p \ dx  \right]^{1/p},
  $$

  $$
  Df = \{ \partial f/\partial x_1, \partial f/\partial x_2, \partial f/\partial x_3,
\ldots,  \partial f/\partial x_m \} = \grad_x f,
  $$

$$
|Df|_p = \left| \left[ \sum_{i=1}^m ( \partial f/\partial x_i)^2  )  \right]^{1/2} \right|_p.
$$
The best possible constant in the inequality (1) belongs to G.Talenti \cite{Talenti1}:

$$
K_m(p) = \pi^{-1/2} m^{-1/p} \left[\frac{p-1}{m-p}\right]^{1-1/p} \cdot
\left[\frac{\Gamma(1+m/2) \ \Gamma(m)}{\Gamma(m/p) \ \Gamma(1+m-m/p)} \right]^{1/m}.
$$

\vspace{3mm}

{\bf B. Trace Sobolev's inequality.} \par

\vspace{3mm}

Let $ m,n = 1,2,\ldots, \ x \in R^m, y \in R^n, z = \{x,y \} \in R^{m+n}, \ u =
u(x,y) = u(z)  $ be any function from the space $ W^1_p(R^{m+n}). $ \par
 We consider in this case only the so-called {\it radial } functions. In detail, we
 define as usually for the vectors $ x = \vec{x} = \{x_1,x_2, \ldots, x_m \} $ and
 $ y = \vec{y} = \{y_1,y_2, \ldots, y_n \} $

 $$
 |x| = \left(\sum_{i=1}^m (x_i)^2 \right)^{1/2}, \ |y| =
 \left(\sum_{j=1}^n (y_j)^2   \right)^{1/2}
 $$
 and correspondingly

 $$
 |z| = |(x,y)| = \left(|x|^2 + |y|^2 \right)^{1/2}.
 $$

We assume that the function $ u = u(x,y) $ depended only on the variable $ |z| $ and
we will write for simplicity

$$
u(x,y) = u(|z|).
$$

Let us denote $ N = m + n, (N \ge 3); \ S[u](x) = u(x,0),$
$$
 \nabla u = \{ \partial u/\partial x_1, \partial u/\partial x_2, \ldots,
 \partial u/\partial x_m, \partial u/\partial y_1, \partial u/\partial  y_2,
  \ldots,  \partial u/\partial y_n \}=
$$

$$
\{\grad_x u, \ \grad_y u \}; \ |\nabla u|_p = (|\grad_x u|^p_p +
 |\grad_y u|^p_p)^{1/p}.
$$

We will denote the class of all the radial functions $ \Rad = \Rad(R^N); \ u(\cdot)
\in \Rad. $ \par

 Notice that the operator $ S[u] $ is correct and continuously defined in the
 $ L_p(R^m) $ in the following sense:

 $$
\lim_{|y| \to 0} |u(\cdot,y) - S[u]|L_p(R^m) = 0,
 $$
see \cite{Besov1},  chapter 5, section 24.\par

The following inequality is called the {\it  Sobolev's trace inequality:}

$$
|S[u](\cdot)|_{q,m} \le K_{m,n}(p) \cdot |\nabla u|_{p,N}, \ q = q(p) = mp/(N-p), \
p \in [1,N). \eqno(2).
$$
   We will understand further under the constant $ K_{m,n}(p) $ in the inequality
(2) its minimal value, namely:

$$
K_{m,n}(p)= \sup
\left\{ \left[ \frac{|S[u](\cdot)|_q}{|\nabla u|_p} \right], \
u \in W^1_p(R^{m+n})\cap \Rad(R^N), \nabla u \ne 0 \right\}. \eqno(3)
$$

 It is evident $ K(m,0) = K(m). $ \par
More information about the constant $ K_{m,n}(p)$ see, for instance, in the articles
\cite{Beesack1}, \cite{Nazaret1}, \cite{Young1}, \cite{Biezuner1}, \cite{Druet1}, \cite{Lieb1}, \cite{Zhu1}, \cite{Edmunds1}, \cite{Escobar1}, \cite{Maz'ja1}
 etc., see also reference therein.\par

\vspace{3mm}

{\bf Our aim is generalization of Sobolev's-type inequalities (1), (3) on some
popular classes of rearrangement invariant (r.i.) spaces, namely, on the so-called
 Grand Lebesgue Spaces $ G(\psi). $  We intend to show also the exactness of offered estimations.}

 \vspace{4mm}

 Hereafter $ C, C_j $ will denote any non-essential finite positive constants.
 We define also for the values $ (p_1, p_2), $ where $ 1 \le p_1 < p_2 \le \infty $

$$
L(p_1, p_2) =  \cap_{p \in (p_1, p_2)} \ L_p.
$$

\vspace{3mm}

  The paper is organized as follows. In the next section we recall the definition
and some simple properties of the so-called Grand Lebesgue Spaces $ G(\psi). $
 In the section 3 we formulate and prove the main result: the classical Sobolev's
inequality  for $ G(\psi) $ spaces  with the exact constant computation. \par
 In the section 4 we investigate the trace Sobolev's inequality for radial functions,
also with the exact constant computation.

  The last section contains some concluding remarks: a {\it weight} generalizations of
ordinary and trace Sobolev's inequality  \par

\vspace{4mm}

\section{Auxiliary facts. Grand Lebesgue Spaces.}\par

 {\bf Definition.} \par

\vspace{3mm}

     Recently, see \cite{Kozachenko1}, \cite{Fiorenza1},
     \cite{Fiorenza2}, \cite{Fiorenza3}, \cite{Iwaniec1}, \cite{Iwaniec2},
     \cite{Ostrovsky1}, \cite{Ostrovsky2}, \cite{Ostrovsky3}, \cite{Ostrovsky4},
     \cite{Ostrovsky5}, \cite{Ostrovsky6}  etc.
     appears the so-called Grand Lebesgue Spaces $ GLS = G(\psi) =
    G(\psi; A,B), \ A,B = \const, A \ge 1, A < B \le \infty, $ spaces consisting
    on all the measurable functions $ f: T \to R $ with finite norms

     $$
     ||f||G(\psi) \stackrel{def}{=} \sup_{p \in (A,B)} \left[ |f|_p /\psi(p) \right].
     $$

      Here $ \psi(\cdot) $ is some continuous positive on the {\it open} interval
    $ (A,B) $ function such that

     $$
     \inf_{p \in (A,B)} \psi(p) > 0, \ \sup_{p \in (A,B)} \psi(p) = \infty.
     $$
We can suppose without loss of generality
$$
\inf_{p \in (A,B)} \psi(p) = 1.
$$

  This spaces are rearrangement invariant, see \cite{Bennet1}, and
  are used, for example, in the theory of probability \cite{Talagrand1}, \cite{Kozachenko1}, \cite{Ostrovsky1}; theory of Partial Differential Equations \cite{Fiorenza2}, \cite{Iwaniec2}; functional analysis \cite{Ostrovsky4}, \cite{Ostrovsky5}; theory of Fourier series \cite{Ostrovsky7}, theory of martingales \cite{Ostrovsky2} etc.\par

\vspace{3mm}

  \section{Classical Sobolev's inequality for Grand Lebesgue Spaces.}

\vspace{3mm}

 We recall $ q = q(p)= mp/(m-p), \ p \in [1,m), $ or equally

$$
p = p(q) = \frac{mq}{m+q}, \ q \in [m/(m-1), \infty).
$$

Let $ \psi(\cdot) $ be arbitrary function from the set $G\Psi(1,m). $ We define the new
function  $ \nu_{\psi}(q) = \nu(q) $ from the set $ G\Psi(m/(m-1), \infty) $ as follows:

$$
\nu_{\psi}(q) = K_m \left( \frac{mq}{m+q} \right) \cdot
\psi\left( \frac{mq}{m+q} \right), \ q \in [m/(m-1), \infty). \eqno(4)
$$

{\bf Theorem 1}. {\it The following Sobolev's type inequality holds:}

$$
||u||G(\nu_{\psi}) \le 1 \cdot ||\nabla u||G(\psi), \eqno(5)
$$
{\it  and the constant "one" in the inequality (5) is the best possible.} \par
{\bf Remark 1.} It is presumed in the assertion (5) of theorem 1 that the right-hand
 side of the proposition of theorem 1 is finite. \par
{\bf Proof.} \par
{\bf A. The upper bound.} Let the function $ u = u(x) $ be such that

$$
\nabla u \in \cap_{p \in (1,m)} W^1_p(R^{m})
$$
and
$$
||\nabla u||G(\psi) \in (0,\infty);
$$
we can assume without loss of generality
$$
||\nabla u||G(\psi) = 1.
$$
It follows from the direct definition of norm for the $ G(\psi) $ spaces that

$$
|\nabla u|_p \le \psi(p), \ p \in (1,m).
$$

We use  now the Talenti's inequality (2):

$$
|u|_{pm/(m-p)} \le K_m(p) \cdot \psi(p), \ p \in (1,m). \eqno(6)
$$
The proposition of theorem 1 follows from (6) after the substitution
$$
p = \frac{mq}{m+q}, \ q \in (m(m-1), \infty):
$$

$$
|u|_{q} \le K_m \left(\frac{mq}{m+q} \right) \cdot \psi
\left(\frac{mq}{m+q}  \right) = \nu(q) = \nu(q) \cdot ||\nabla u||G(\psi);
$$

$$
\frac{|u|_{q}}{\nu(q)} \le ||\nabla u||G(\psi);
$$
therefore

$$
||u||G(\nu_{\psi}) = \sup_q \left[\frac{|u|_{q}}{\nu(q)} \right]
 \le ||\nabla u||G(\psi).
$$

\vspace{3mm}

{\bf B. Proof of the low bound.}\par
{\bf 1.} Let us introduce the following important functional: $ V  = $

$$
\sup_{m=3,4,\ldots} \sup_{\psi \in G\Psi(m/(m-1),m)} \sup
  \left\{ \frac{||f||G(\nu)}{||\nabla f||G(\psi)}: f \in G(\psi)\cap W^1_p(R^m),
  \nabla f \ne 0 \right\}. \eqno(7)
$$
We proved that $ V \le 1; $ it remains to ground the opposite inequality.\par

 Let us consider the following example (more exactly, a family of examples)
 of a functions:

$$
 f_{\Delta} = f_{\Delta}(x) = |\log |x| \ |^{\Delta} \cdot I(|x| \le 1),
$$
where $ \Delta = \const \ge 1, $

$$
I(|x| \le 1) = 1, \ |x| \le 1; I(|x| \le 1) = 0, \ |x| >1.
$$

So, $ f_{\Delta}(x) $ is radial function and  $ f_{\Delta}(\cdot) \in
G(\psi)\cap W^1_p(R^m), \nabla f_{\Delta}(\cdot) \ne 0. $ \par
 We calculate denoting

 $$
 \omega(m)= \frac{2 \pi^{m/2}}{\Gamma(m/2)}:
 $$

 $$
 |f_{\Delta}|_q =
 \frac{\omega^{1/q}(m) \ \Gamma^{1/q}(\Delta q + 1)}{m^{\Delta + 1/q}};
 $$

$$
|\nabla f_{\Delta}|_p = \Delta \cdot \omega^{1/p}(m) \cdot
\frac{\Gamma^{1/p}(p(\Delta-1)+1)}{(m-p)^{\Delta-1+1/p}}.
$$
 We can rewrite the expression for the value $ V $ as follows: $ V = $

$$
\sup_m \sup_{\psi \in G\Psi(m/(m-1),m)} \sup
  \left\{ \frac{ \sup_q[|f|_q/\nu(q)]}{\sup_p [|\nabla f|_p/\psi(p)]}: f \in G(\psi)\cap W^1_p(R^m),  \nabla f \ne 0 \right\}. \eqno(8)
$$
 When we choose $ f = f_{\Delta} $ and $ \psi(p) = \psi_{\Delta}(p) :=
 | \nabla f_{\Delta}|_p, $ we obtain the following {\it low} bond for the value
 $ V: V \ge V_0,$  where

$$
V_0 = \sup_m \sup_{p \in (1,m)} \left[ \frac{ |f_{\Delta}|_{q(p)}}{K_m(p) \ |\nabla f_{\Delta}|_p } \right]. \eqno(9)
$$

 Substituting the expressions for $ |f_{\Delta}|_q, \ |\nabla f_{\Delta}|_p $ and for
 $ K_m(p) $ into the formula (9), we obtain after some calculations by means of
 Stirling's formula for all the admissible values $m, \Delta: $

 $$
 V_0 \ge \lim_{p \to m-0} \left[ \frac{ |f_{\Delta}|_{q(p)}}{K_m(p) \ |\nabla f_{\Delta}|_p } \right] =
 \frac{m^{\Delta} \Delta^{\Delta-1} e^{-\Delta}  \ m^{1/m}}
 { (m-1)^{1-1/m} \ \Gamma^{1/m}((\Delta-1)m + 1)}=:V_{00}(m,\Delta);
 $$

 $$
 V_0 \ge \lim_{m \to \infty} V_{00}(m,\Delta) = e^{-1}
  \left[\frac{\Delta}{\Delta -1} \right]^{\Delta-1} =:V_{000}(\Delta).
 $$
 Finally,

 $$
 V_0 \ge \lim_{\Delta \to \infty} V_{000}(\Delta) = 1.
 $$
 This completes the proof of theorem 1.\par

\vspace{3mm}

\section{Trace Sobolev's inequality for Grand Lebesgue Spaces. \ Radial case.}

\vspace{3mm}

We recall that we consider in this section {\it only radial } functions
$ u(z) = u(x,y) = u(|z|) $ and that here

$$
q = q(p)= mp/(N-p), \ p \in [p_0, N], \ p_0 \stackrel{def}{=}
\max(N/(m+1),1),
$$
or equally

$$
p = p(q) = \frac{mq}{N+q}, \ q \in [1, \infty).
$$

Let $ \psi(\cdot) $ be arbitrary function from the set $G\Psi(p_0,N). $
We define the new function  $ \zeta_{\psi}(q) = \zeta(q) $ from the set
$ G\Psi(1, \infty) $ as follows:

$$
\zeta_{\psi}(q) = K_m \left( \frac{mq}{m+q} \right) \cdot
\psi\left( \frac{mq}{m+q} \right), \ q \in [1, \infty). \eqno(10)
$$

{\bf Theorem 2}. {\it The following trace Sobolev's type inequality holds:}

$$
||u||G(\nu_{\psi}) \le 1 \cdot ||\nabla u||G(\psi), \eqno(11)
$$
{\it  and the constant "one" in the inequality (11) is the best possible.} \par

{\bf Remark 2.} As in the remark 1,
it is presumed in the assertion (11) of theorem 2 that the right-hand
 side of the proposition of theorem 2 is finite. \par
{\bf Proof } is alike to the proof of the assertion of theorem 1. \par
{\bf A. The upper bound.} Let the function $ u = u(z) = u(|z|) $ be radial function:
$ u(\cdot) \in Rad(R^N) $ and such that

$$
\nabla u \in \cap_{p \in (p_0,N)} W^1_p(R^{N}), \ \nabla u \ne 0.
$$

We assume without loss of generality that
$$
||\nabla u||G(\psi) = 1.
$$

It follows from the definition of norm for the $ G(\psi) $ spaces that
$$
|\nabla u|_p \le \psi(p), \ p \in [p_0,N).
$$

It follows from the definition of the constants $ K_{m,n}(p) $
(inequality (3)):

$$
|S[u]|_{pm/(N-p)} \le K_{m,n}(p) \cdot \psi(p), \ p \in (p_0,N). \eqno(12)
$$
The proposition of theorem 2 follows from (12) after the substitution
$$
p = \frac{Nq}{m+q}, \ q \in (1, \infty):
$$

$$
|S[u]|_{q} \le K_{m,n} \left(\frac{Nq}{m+q} \right) \cdot \psi
\left(\frac{Nq}{m+q}  \right) = \zeta(q) = \zeta(q) \cdot ||\nabla u||G(\psi);
$$

$$
\frac{S[u]_{q}}{\zeta(q)} \le ||\nabla u||G(\psi);
$$
thus

$$
||S[u]||G(\zeta_{\psi}) = \sup_q \left[\frac{|S[u]|_{q}}{\zeta(q)} \right]
 \le ||\nabla u||G(\psi).
$$

\vspace{3mm}

{\bf B. Proof of the low bound.}\par
\vspace{3mm}
{\bf 1.} Let us introduce again the following important functional
(with at the same notation as in the last section): $ V(m,n) = $

$$
\sup_N \sup_{\psi \in G\Psi(p_0,N)} \sup
  \left\{ \frac{||S[f]||G(\zeta)}{||\nabla f||G(\psi)}: \ \nabla f \in G(\psi)\cap W^1_p(R^N),  \nabla f \ne 0 \right\}. \eqno(13)
$$
We proved that $ V(m,n) \le 1; $ it remains to ground an opposite inequality.\par

\vspace{3mm}

{\bf 2.} Let us consider again the used examples of a functions:
$$
 f_{\Delta} = f_{\Delta}(z) = |\log |z| \ |^{\Delta} \cdot I(|z| \le 1), \eqno(14)
$$
where $ \Delta = \const \ge 1. $ \par

So, $ f_{\Delta}(x) $ is radial function and such that  $ \nabla f_{\Delta}(\cdot) \in
G(\psi)\cap W^1_p(R^N), \nabla f_{\Delta}(\cdot) \ne 0. $ \par
 We calculate as before:

  $$
 |S[f_{\Delta}]|_q =
 \frac{\omega^{1/q}(m) \ \Gamma^{1/q}(\Delta q + 1)}{m^{\Delta + 1/q}}; \eqno(15)
 $$

$$
|\nabla f_{\Delta}|_p = \Delta \cdot \omega^{1/p}(N) \cdot
\frac{\Gamma^{1/p}(p(\Delta-1)+1)}{(N-p)^{\Delta-1+1/p}}. \eqno(16)
$$

\vspace{3mm}

{\bf 3.} We do not know, in contradiction to the case of ordinary Sobolev's inequality,
the exact value of a constant $ K_{m,n}(p). $  In order to prove the proposition of
theorem 2, we need to obtain the {\it upper estimate } for this constant. \par

 Let us estimate the constant $ K_{m,n}(p). $ As long as the function
$ f(\cdot)$ is radial, we can rewrite the inequality (2) using the multidimensional
spherical coordinates as follows:

$$
\omega^{1/q}(m) \left[ \int_0^{\infty} \left|s^{m-1} \ ds \cdot
 \int_s^{\infty} g(t) \ dt \right| \  \right]^{1/q} \le \omega^{1/p}(N) \ K_{m,n}(p) \times
$$

$$
\left[\int_0^{\infty} s^{N-1} |g(s)|^p \ ds \right]^{1/p}. \eqno(17)
$$
 Further we will use the result belonging to Bradley \cite{Bradley1};
 see also \cite{Maz'ja1}, which is a weight generalization of the classical
 Hardy-Littlewood inequality. It asserts that the following inequality is true:

 $$
 \left\{ \left[ \int_0^{\infty} u(x) \int_x^{\infty} f(t)dt \right]^q \ dx
 \right\}^{1/q} \le
C \times \left\{ \int_0^{\infty} [v(x) f(x)]^p \ dx \right\}^{1/p},  \eqno(18)
$$
where $ u(x), v(x) \ge 0, 1 < p \le q < \infty,  $

$$
Q(p):= p^{1/q} \cdot ( p/(p-1))^{(p-1)/p },
$$

$$
B \le C \le B \cdot Q(p),
$$

$$
B = \sup_{w>0}J(w), \ J(w) = J_1(w) \cdot J_2(w),
$$

$$
J_1(w) =  \left( \int_0^w u^q(x) dx \right)^{1/q};
J_2(w) = \left( \int_w^{\infty} (v(x))^{-(p-1)/p} \ dx  \right)^{(p-1)/p}. \eqno(19)
$$
 Note that the case $ p \le q, $ i.e. $ p \ge n $ is sufficient in order to prove the
 second assertion of theorem 2, as long as we put further $ p \to N-0. $ \par

 We compute in the considered case:

$$
J_1(w) = m^{-1/q} \ w^{m/q}, \ J_2(w) = \left[ \frac{p-1}{N-p} \right]^{1-1/p} \cdot
w^{-(N-p)/p}.
$$
Hence, the expression for the value of $ B $ is finite only in the case when

$$
\frac{m}{q}= \frac{N-p}{p},
$$
or equally
$$
q = \frac{mp}{N-p},
$$
i.e. as in the conditions of theorem 2. \par
We conclude in the considered case:

$$
B = m^{-1/q} \ \left[\frac{p-1}{N-p} \right]^{1-1/p}
$$
and following

$$
 C \le Q(p) \cdot m^{-1/q} \cdot \left[ \frac{p-1}{N-p} \right]^{1-1/p},
$$

$$
K_{m,n}(p) \le K^+_{m,n}(p), \ K^+_{m,n}(p) \stackrel{def}{=} \omega^{1/q}(m) \cdot
 \omega^{-1/p}(N) \cdot C =
$$

$$
Q(p) \cdot m^{-1/q} \cdot \left[\frac{p-1}{N-p} \right]^{1-1/p} \cdot
\omega^{1/q}(m) \cdot \omega^{-1/p}(N). \eqno(20)
$$

\vspace{3mm}

{\bf 4.} We can rewrite the expression for the value $ V(m,n) $ as follows:
$ V(m,n) = $

$$
\sup_N  \sup_{\psi \in G\Psi(m/(m-1),m)} \sup
 \left\{ \frac{ \sup_q[|f|_q/\zeta(q)]}{\sup_p [|\nabla f|_p/\psi(p)]}:
 \nabla f \in G(\psi)\cap W^1_p(R^m),  \nabla f \ne 0 \right\}. \eqno(21)
$$
 When we choose $ f = f_{\Delta} $ and $ \psi(p) = \psi_{\Delta}(p) :=
 | \nabla f_{\Delta}|_p, $ we obtain the following {\it low} bond for the value
 $ V(m,n): V(m,n) \ge V_0(m,n),$  where

$$
V_0(m,n) =\sup_{\Delta \in (1,\infty)} \sup_N  \sup_{p \in (p_0,m)}
\left[ \frac{ |f_{\Delta}|_{q(p)}}{K_{m,n}(p) \ |\nabla f_{\Delta}|_p } \right].
$$

 Substituting the expressions for $ |f_{\Delta}|_q, \ |\nabla f_{\Delta}|_p $
 into the formula for the value $ V(m,n),$ tacking into account the
 inequality $ K_{m,n}(p) \le K^+_{m,n}(p), $ we obtain after some calculations by
 means of Stirling's formula for all the admissible values $ m,p, \Delta: $

 $$
 V_0(m,n) \ge \lim_{p \to N-0} \left[ \frac{ |f_{\Delta}|_{q(p)}}{K^+_{m,n}(p)
  \ |\nabla f_{\Delta}|_p } \right] =
$$

$$
 \frac{\Delta^{\Delta-1} e^{-\Delta}  N^{\Delta-1}}
 {\Gamma^{1/N}((\Delta-1)N + 1)}=:V_{00}(m,n;\Delta,N); \eqno(22)
 $$

 $$
 V_0(m,n) \ge \lim_{N \to \infty} V_{00}(m,n;\Delta,N) = e^{-1}
  \left[\frac{\Delta}{\Delta -1} \right]^{\Delta-1} =:V_{000}(m,n;\Delta). \eqno(23)
 $$
 Finally,
$$
 V_0(m,n) \ge \lim_{\Delta \to \infty} V_{000}(m,n)(\Delta) = 1. \eqno(24)
 $$
 This completes the proof of theorem 2.\par

  Note that at the same result may be obtained from the estimation in \cite{Besov2},
 chapter 1, section 2:

 $$
 \omega^{1/q}(m) |S[u]|_q  \le \omega^{1/p}(N) \cdot
 \left[\frac{p(m-1) + N}{m(N-p)} \right]^{1-(p-n)/(mp)}\cdot |\nabla u|_p, \eqno(25)
 $$
hence
$$
K_{m,n}(p) \le K^{(1)}_{m,n}(p) \stackrel{def}{=} \omega^{1/p}(N) \cdot \omega^{-1/q}(m) \cdot \left[\frac{p(m-1)+ N}{m(N-p)} \right]^{1-(p-n)/(mp)}, \eqno(26)
$$

$$
V(m,n) \ge \lim_{\Delta \to \infty} \lim_{N \to \infty} \lim_{p \to N-0}
\left[\frac{|S[f_{\Delta}]_{q(p)}}{K^{(1)}_{m,n} (p) \cdot
 |\nabla f_{\Delta}|_p } \right] = 1.
$$

\vspace{3mm}

\section{ Concluding remarks. Weight Sobolev's inequalities.}

\vspace{3mm}

{\bf A.} We introduce the so-called {\it weight trace operator} by the formula 

$$
 S_{\alpha}[u](x) =  |x|^{-\alpha} \cdot u(x,0), \ u(z) = u(x,y) = u(|z|), 
$$
i.e. $ u(\cdot) \in Rad(R^N). $ \par
 Let us consider a {\it weight} Poincare-Sobolev trace inequality of a view:

$$
| S_{\alpha} [u]|_q \le K^{\alpha}_{m,n}(p) \cdot |\nabla u|_p, \eqno(27)
$$
where
$$
\alpha = \const \in [0,1], \ p \in (p_1, N/(1-\alpha)),
$$

$$
p_1 \stackrel{def}{=} \max(N/(m+1-\alpha), 1) < N/(1-\alpha); \ N/0 = + \infty;
\eqno(28)
$$

$$
q = q(p) = \frac{mp}{N-p(1-\alpha)}; \ q \in (1,\infty) \Leftrightarrow
$$

$$
p = \frac{qN}{m+q(1-\alpha)}. \eqno(29)
$$

and we understood the value $ K^{\alpha}_{m,n}(p) $ as the its minimal value:
$$
K^{\alpha}_{m,n}(p) = \sup_{u: \nabla u \ne 0}
\left[\frac{| S_{-\alpha}[u]|_q}{|\nabla u|_p}\right] < \infty, 
p \in (p_1,N/(1-\alpha)).
$$
 Notice that in the case $ \alpha = 0 $ we obtain the Sobolev's inequality and
 the case $ \alpha = 1 $ correspondent the so-called modified Poincare's inequality.
 \par

Let $ \psi \in G\Psi(p_1, N/(1-\alpha)); $ we introduce a new function

$$
\theta_{\psi}(q) = \psi \left(\frac{qN}{m+q(1-\alpha)} \right) \cdot
K^{\alpha}_{m,n} \left(\frac{qN}{m+q(1-\alpha)} \right). \eqno(30)
$$

{\bf Theorem 3}. {\it The following generalized trace Sobolev-Poincare type inequality holds:}

$$
||S_{\alpha}[u] ||G(\theta_{\psi}) \le 1 \cdot ||\nabla u||G(\psi), \eqno(31)
$$
{\it  and the constant "one" in the last inequality  is the best possible.} \par

{\bf Proof} is at the same as in the theorem 2, with at the same "counterexamples"
and may be omitted. \par

\vspace{3mm}

{\bf B.} Note that the {\it low bound} for Sobolev's trace embedding constants, for
instance, $ V(m,n) \ge 1, $ are true for arbitrary, i.e. not only  for radial functions. \par

\vspace{3mm}
{\bf C.} The exact value for the degree $ q $ in the generalized Poincare-Sobolev inequality, namely

$$
q = q(p) = \frac{mp}{N-p(1-\alpha)}
$$
may be obtained by means of the so-called {\it dilation method}, offered by Talenti
\cite{Talenti1}. In detail, let us define as usually the family of dilation operators
$ T_{\lambda}[f], \ \lambda \in (0,\infty), \ f: R^k \to R, \ k = 1,2,\ldots, $
of a view:

$$
T_{\lambda}[f](x) = f(x/\lambda). \eqno(32)
$$
Suppose the inequality (27) is satisfied for {\it some} admissible radial function such that $ \nabla u \in W^1_p(R^N), \nabla u \ne 0.$ 
As long as $ u \in W^1_p(R^N) \Rightarrow T_{\lambda}u \in
 W^1_p(R^N), $ we  have rewriting the inequality (27) for the function $ T_{\lambda}u: $
 
$$
| S_{\alpha} [T_{\lambda}u]|_q \le K^{\alpha}_{m,n}(p) \cdot 
|\nabla T_{\lambda}u|_p, \eqno(33)
$$
tacking into account the equalities:

$$
|T_{\lambda}u|_{p,N} = \lambda^{N/p} |u|_{p,N},
$$

$$
|S_{\alpha} [T_{\lambda}u]|_{q,m} = \lambda^{m/q - \alpha} |S[u]|_q,
$$

$$
\nabla T_{\lambda}[ u](z) = \lambda^{-1} u(z/\lambda), 
$$

$$
|\nabla T_{\lambda}[ u]|_{p,N} = \lambda^{N/p -1} \ |u|_{p,N}:
$$

$$
\lambda^{m/q - \alpha} |S_{\alpha}[u]|_{q,m} \le K^{\alpha}_{m,n}(p) \ \lambda^{N/p-1} \ 
|\nabla u|_{p,N}. \eqno(34)
$$
 Since the value $ \lambda $ in the last inequality  is arbitrary in the set 
 $ (0,\infty), $ we conclude
 $$
 m/q - \alpha = N/p - 1
 $$
or equally 

$$
q = \frac{p \ m}{N-p(1-\alpha)}.
$$


\vspace{3mm}

\begin{thebibliography}{99}

\vspace{3mm}

\bibitem{Beesack1}
{\sc Beesack P.R.} Hardy's inequality and its extension. Pacific J. Math.;
{\bf 11}, (1961), 39-61.
\bibitem{Bennet1}
 {\sc Bennet G., Sharpley R.} Interpolation of operators. Orlando, Academic Press Inc., (1988).
\bibitem{Besov1}
{\sc O.V.Besov, P.I.Il'in, S.M. Nikol'skii.} Integral Representations of Functions
and Imbedding Theorems. Volume 1, (1978), John Wiley and Sons, Washington D.C., New York, Toronto, London, Sydney; A Halisted Press Book.
\bibitem{Besov2}
{\sc O.V.Besov, P.I.Il'in, S.M. Nikol'skii.} Integral Representations of Functions
and Imbedding Theorems. Volume 2, (1979), John Wiley and Sons, Washington D.C., New York, Toronto, London, Sydney; A Halisted Press Book.
\bibitem{Biezuner1}
{\sc R.J.Biezuner.} Best constants in Sobolev trace inequalities. Nonlinear Analysis,
{\bf 54}, (2003), 457-502.
\bibitem{Bradley1}
{\sc J.S.Bradley.}  Hardy inequalities with mixed norms. Canadian Math. Bull., 21(1978), p. 405-408.
\bibitem{Druet1}
{\sc O.Druet.} The best constants problem in Sobolev trace inequalities. Math. Annalen,
{\bf 314}, (1999), 327-346.
\bibitem{Edmunds1}
{\sc D.E.Edmunds and W.D. Evans.}  Sobolev Embeddings and Hardy Operators. In:
 Vladimir Maz'ya (Editor), "Sobolev Spaces in Mathematics", Part 1, International
 Mathematical Series,  Volume 8, Springer Verlag, Tamara Rozhkovskaya Publisher;
 (2009), New York, London, Berlin; p.153-184.
\bibitem{Escobar1}
{\sc J.F.Escobar.} Sharp constant in a Sobolev trace inequality. Indiana Math, J.,
{\bf 37}, (1988), 687-698.
\bibitem{Fiorenza1}
   {\sc Capone C., Fiorenza A., Krbec M.} On the Extrapolation Blowups in the
   $ L_p $ Scale. Collectanea Mathematica, {\bf 48}, 2, (1998), 71 - 88.
 \bibitem{Fiorenza2}
 {\sc A.Fiorenza.} Duality and reflexivity in grand Lebesgue spaces.
       Collectanea Mathematica (electronic version), {\bf 51}, 2, (2000), 131 - 148.
\bibitem{Fiorenza3}
 {\sc A. Fiorenza and G.E. Karadzhov.} Grand and small Lebesgue spaces and
       their analogs. Consiglio Nationale Delle Ricerche, Instituto per le
      Applicazioni del Calcoto Mauro Picine", Sezione di Napoli, Rapporto tecnico n.
      272/03, (2005).
\bibitem{Iwaniec1}
   {\sc T.Iwaniec and C. Sbordone.} On the integrability of the Jacobian under
      minimal hypotheses. Arch. Rat.Mech. Anal., 119, (1992), 129 – 143.
\bibitem{Iwaniec2}
 {\sc T.Iwaniec, P. Koskela and J. Onninen.} Mapping of finite distortion:
   Monotonicity and Continuity.  Invent. Math. 144 (2001), 507 - 531.
\bibitem{Kantorovicz1}
{\sc L.V.Kantorovicz, G.P.Akilov.} Functional Analysis. (1987) Kluvner Verlag.
\bibitem{Kozachenko1}
 {\sc Kozachenko Yu. V., Ostrovsky E.I.} (1985). The Banach Spaces of
      random Variables of subgaussian type. {\it Theory of Probab. and Math.
      Stat.} (in Russian). Kiev, KSU, {\bf 32}, 43 - 57.
\bibitem{Talagrand1}
 {\sc Ledoux M., Talagrand M.} (1991) Probability in Banach Spaces.
  Springer, Berlin, MR 1102015.
\bibitem{Lieb1}
{\sc E.H.Lieb.} Sharp constants in the Hardy-Littlewood-Sobolev and related inequalities.
Ann. Math., {\it 118},  (1983), 349-374.
\bibitem{Maz'ja1}
{\sc V.Maz'ja.} Sobolev Spaces. Kluvner Academic Verlag, (2002),
Berlin-Heidelberg-New York.
\bibitem{Mitrinovich1}
{\sc D.S.Mitrinovich, J.E. Pecaric and A.M.Fink.} Inequalities Involving
Functions and Their Integrals and Derivatives. Kluvner Academic Verlag,
(1996), Dorderecht, Boston, London.
\bibitem{Nazaret1}
{\sc B.Nazaret.} Best constant in Sobolev trace inequalities on the half space.
arXiv:math.FA/2109077 v1 21Aug 2005.
\bibitem{Ostrovsky1}
{\sc E.I. Ostrovsky.}  Exponential Estimations for Random Fields.
Moscow - Obninsk, OINPE, 1999 (Russian).
\bibitem{Ostrovsky2}
{\sc E. Ostrovsky and L.Sirota.} Moment Banach spaces: theory and applications.
HAIT Journal of Science and Engineering, {\bf C}, Volume 4, Issues 1 - 2,
pp. 233 - 262, (2007).
\bibitem{Ostrovsky3}
{\sc E.Ostrovsky, E.Rogover and L.Sirota.} Riesz's and Bessel's operators in in
bilateral Grand Lebesgue Spaces.
arXiv:0907.3321 [math.FA] 19 Jul 2009.
\bibitem{Ostrovsky4}
{\sc E. Ostrovsky and L.Sirota.} Weight Hardy-Littlewood inequalities for different
powers. arXiv:09010.4609v1[math.FA] 29 Oct 2009.
\bibitem{Ostrovsky5}
{\sc E. Ostrovsky E.} Bide-side exponential and moment inequalities
     for tail of distribution of Polynomial Martingales. Electronic
     publication, arXiv: math.PR/0406532 v.1 Jun. 2004.
\bibitem{Ostrovsky6}
{\sc E.Ostrovsky,  E.Rogover and L.Sirota.} Integral Operators in Bilateral Grand
Lebesgue Spaces. arXiv:09012.7601v1 [math.FA] 16 Dez. 2009.
\bibitem{Ostrovsky7}
{\sc E.Ostrovsky, L.Sirota.} Nikolskii-type inequalities for rearrangement invariant
spaces. arXiv:0804.2311v1 [math.FA] 15 Apr 2008.
\bibitem{Ostrovsky8}
 {\sc Ostrovsky E., Sirota L.} Moment Banach Spaces: Theory and Applications.
 HIAT Journal of Science and Engineering, Holon, Israel, v. 4, Issue 1-2, (2007), 233 - 262.
\bibitem{Sobolev1}
{\sc Sobolev S.L.} Some Applications of Functional Analysis into Mathematical Physic,
(1950), Publishing House LSU (Leningrad State University), (in Russian).
\bibitem{Talenti1}
{\sc G.Talenti} Inequalities in Rearrangement Invariant Function Spaces. Nonlinear Analysis, Function Spaces and Applications. Prometheus, Prague, {\bf 5}, (1995), 177-230.
\bibitem{Young1}
{\sc Young Ja Park.} Sobolev Trace Inequalities.
arXiv:math.CA/0107065 v1 9Jul 2001.
\bibitem{Zhu1}.
{\sc M.Zhu} Some general forms of sharp Sobolev inequalities. J. Func. Anal., {\bf 156}, (1998), 75-120.

\end{thebibliography}
\end{document}